\documentstyle[amscd,amssymb,verbatim,epsf]{amsart}
\begin{document}

\theoremstyle{plain}
\newtheorem{Thm}{Theorem}
\newtheorem{Cor}{Corollary}
\newtheorem{Con}{Conjecture}
\newtheorem{Main}{Main Theorem}

\newtheorem{Lem}{Lemma}
\newtheorem{Prop}{Proposition}

\theoremstyle{definition}
\newtheorem{Def}{Definition}
\newtheorem{Note}{Note}

\theoremstyle{remark}
\newtheorem{notation}{Notation}
\renewcommand{\thenotation}{}

\errorcontextlines=0
\numberwithin{equation}{section}
\renewcommand{\rm}{\normalshape}%

\title[Half quasi-Cauchy sequences]%
   {Half quasi-Cauchy sequences}
\author{H\"Usey\.{I}n \c{C}akall\i\\
          Department of Mathematics, Maltepe University, Marmara E\u{g}\.{I}t\.{I}m K\"oy\"u, TR 34857, Maltepe, \.{I}stanbul-Turkey \; \; \; \; \; Phone:(+90216)6261050 ext:2248, fax:(+90216)6261113}
\address{H\"Usey\.{I}n \c{C}akall\i\\
          Department of Mathematics, Maltepe University, Marmara E\u{g}\.{I}t\.{I}m K\"oy\"u, TR 34857, Maltepe, \.{I}stanbul-Turkey \; \; \; \; \; Phone:(+90216)6261050 ext:2248, fax:(+90216)6261113
}
\email{hcakalli@@maltepe.edu.tr}

\date{\today}

\begin{abstract}

A real function $f$ is ward continuous if $f$ preserves quasi-Cauchyness, i.e. $(f(x_{n}))$ is a quasi-Cauchy sequence whenever $(x_{n})$ is quasi-Cauchy; and a subset $E$ of $\textbf{R}$ is quasi-Cauchy compact if any sequence $\textbf{x}=(x_{n})$ of points in $E$ has a quasi-Cauchy subsequence where $\textbf{R}$ is the set of real numbers. These known results suggest to us introducing a concept of upward (respectively, downward) half quasi-Cauchy continuity in the sense that a function $f$ is upward (respectively, downward) half quasi-Cauchy continuous if it preserves upward (respectively, downward) half quasi-Cauchy sequences, and a concept of upward (respectively, downward) half quasi-Cauchy compactness in the sense that a subset $E$ of $\textbf{R}$ is upward (respectively, downward) half quasi-Cauchy compact if any sequence of points in $E$ has an upward (respectively, downward) half quasi-Cauchy subsequence. We investigate upward(respectively, downward) half quasi-Cauchy continuity and upward (respectively, downward) half quasi-Cauchy compactness, and prove related theorems.

\end{abstract}

\maketitle

\section{Introduction}
\normalfont

Recently, a concept of forward continuity of a real function $f$ was introduced in the sense that $f$ is called forward continuous if $(f(x_{n}))$ is a quasi-Cauchy sequence whenever $(x_{n})$ is quasi-Cauchy (see \cite{CakalliForwardcompactness}, and \cite{CakalliForwardcontinuity}). A subset $E$ of $\textbf{R}$, the set of real numbers, is bounded if and only if any sequence of points in $E$ has a quasi-Cauchy subsequence where a sequence $\textbf{x}=(x_{n})$ is quasi-Cauchy if for every $\varepsilon>0$ there exists an $n_{0}\in{\textbf{N}}$ such that $|x_{n}-x_{n+1}| <\varepsilon$ for $n \geq n_0$ where $\textbf{N}$ denotes the set of positive integers (see \cite{BurtonColeman}, \cite{CakalliStatisticalquasiCauchysequences}, and \cite{CakalliStatisticalwardcontinuity}). Using the idea of forward continuity of a real function and the idea of forward compactness, we introduce a concept of upward half quasi-Cauchy continuity in the sense that a function $f$ is upward half quasi-Cauchy continuous if it transforms upward half quasi-Cauchy sequences to upward half quasi-Cauchy sequences, i.e. $(f(x_{n}))$ is upward half quasi-Cauchy whenever $(x_{n})$ is upward half quasi-Cauchy, and a concept of upward half quasi-Cauchy compactness in the sense that a subset $E$ of $\textbf{R}$ is upward half quasi-Cauchy compact if any sequence of points in $E$ has an upward half quasi-Cauchy subsequence where a sequence $\textbf{x}=(x_{n})$ is called upward half quasi-Cauchy if for every $\varepsilon>0$ there exists an $n_{0}\in{\textbf{N}}$ so that $x_{n}-x_{n+1} <\varepsilon$ for $n \geq n_0$. We then similarly state a definition of a concept of downward half quasi-Cauchy continuity in the sense that a function $f$ is downward half quasi-Cauchy continuous if it transforms downward half quasi-Cauchy sequences to downward half quasi-Cauchy sequences, i.e. $(f(x_{n}))$ is downward half quasi-Cauchy whenever $(x_{n})$ is downward half quasi-Cauchy, and a definition of a concept of downward half quasi-Cauchy compactness in the sense that a subset $E$ of $\textbf{R}$ is downward half quasi-Cauchy compact if any sequence of points in $E$ has a downward half quasi-Cauchy subsequence where a sequence $\textbf{x}=(x_{n})$ is called downward half quasi-Cauchy if for every $\varepsilon>0$ there exists an $n_{0}\in{\textbf{N}}$ so that $x_{n+1}-x_{n} <\varepsilon$ for $n \geq n_0$.

The purpose of this paper is to introduce the concepts of upward and downward half quasi-Cauchy continuity of a function and the concepts of upward and downward half quasi-Cauchy compactness of a subset of $\textbf{R}$ and prove related theorems.

\maketitle

\section{Preliminaries}

\normalfont{}

We will use boldface letters $\boldsymbol{\alpha}$, $\bf{x}$, $\bf{y}$, $\bf{z}$, ... for sequences $\boldsymbol{\alpha}=(\alpha_{k})$, $\textbf{x}=(x_{n})$, $\textbf{y}=(y_{n})$, $\textbf{z}=(z_{n})$, ... of points in \textbf{R} for the sake of abbreviation. $s$ and $c$ will denote the set of all sequences, and the set of convergent sequences of points in \textbf{R}.

 A subset of \textbf{R} is compact if and only if it is closed and bounded. A subset $A$ of  \textbf{R} is bounded if $|a|\leq{M}$ for all $a \in{A}$  where  $M$ is a positive  real constant number. This is equivalent to the statement that any sequence of points in $A$ has a Cauchy subsequence. The concept of a Cauchy sequence involves far more than that the distance between successive terms is tending to zero. Nevertheless, sequences which satisfy this weaker property are interesting in their own right.  A sequence $(\alpha _{n})$ of points in \textbf{R} is quasi-Cauchy if $(\Delta \alpha _{n})$ is a null sequence where $\Delta \alpha _{n}=\alpha _{n+1}-\alpha _{n}$. These sequences were named as quasi-Cauchy by Burton and Coleman \cite[page 328]{BurtonColeman}, while they  were called as forward convergent to $0$ sequences in \cite[page 226]{CakalliForwardcontinuity}.

It is known that a sequence $(\alpha _{n})$ of points in \textbf{R} is slowly oscillating if
$$
\lim_{\lambda \rightarrow 1^{+}}\overline{\lim}_{n}\max _{n+1\leq
k\leq [\lambda n]} |
  \alpha _{k}  -\alpha _{n} | =0
$$ where $[\lambda n]$ denotes the integer part of $\lambda n$ (see \cite[Definition 2 page 947]{FDikMDikandCanak}).
Any Cauchy sequence is slowly oscillating, and any slowly oscillating sequence is quasi-Cauchy. There are quasi-Cauchy sequences which are not Cauchy. For example, the sequence  $(\sqrt{n})$ is quasi-Cauchy, but not Cauchy. Any subsequence of a Cauchy sequence is Cauchy. The analogous property fails for quasi-Cauchy sequences, and fails for slowly oscillating sequences as well. A counterexample for the case, quasi-Cauchy, is again the sequence $(a_n)=(\sqrt{n})$ with the subsequence $(a_{n^{2}})=(n)$. A counterexample for the case slowly oscillating is the sequence $(log_{10} n)$ with the subsequence (n). Furthermore we give more examples without neglecting: the sequences $(\sum^{\infty}_{k=1}\frac{1}{n})$, $(ln\;n)$, $(ln\;(ln\;n))$, $(ln\;( ln\;( ln n)))$,..., $(ln\;( ln\;( ln (...( ln n)...))$   and combinations like that are all slowly oscillating, but not Cauchy. The bounded sequence $(cos (6 log(n + 1)))$ is slowly oscillating, but not Cauchy. The sequences $(cos(\pi \sqrt{n}))$ and $(\sum^{k=n}_{k=1}(\frac{1}{k})(\sum^{j=k}_{j=1}\frac{1}{j}))$ are quasi-Cauchy, but slowly oscillating(see \cite{CakalliSlowlyoscillatingcontinuity} and \cite{Vallin}).

By a method of sequential convergence, or briefly a method, we mean a linear function $G$ defined on a subspace of $s$, denoted by $c_{G}$, into \textbf{R}. A sequence $\textbf{x}=(x_{n})$ is said to be $G$-convergent to $\ell$ if $\boldsymbol{\alpha}\in c_{G}$ and $G(\boldsymbol{\alpha})=\ell$ (\cite{ConnorandGrosseSequentialdefinitionsofcontinuityforrealfunctions}). In particular, $\lim$ denotes the limit function $\lim \boldsymbol{\alpha}=\lim_{n}\alpha_{n}$ on the space $c$ of convergent sequences of points in \textbf{R}. A method $G$ is called regular if $c\subset{c_{G}}$, i.e. every convergent sequence $\boldsymbol{\alpha}=(\alpha_{n})$ is $G$-convergent with $G(\boldsymbol{\alpha})=\lim \boldsymbol{\alpha}$.

Consider an infinite matrix $\verb"A"=(a_{nk})^{\infty}_{n,k=1}$ of real numbers. Then, for any sequence $\textbf{x}=(x_{n})$ the sequence $\verb"A"\textbf{x}$ is defined as

\[
\verb"A"\textbf{x}=(\sum^{\infty}_{k=1}a_{nk}x_{k})_{n}
\]
provided that each of the series converges. A sequence $\textbf{x}$ is called $\verb"A"$-convergent (or $\verb"A"$-summable) to $\ell$ if $\verb"A"\textbf{x}$ exists and is convergent with

\[
\lim \verb"A"\textbf{x}=\lim_{n\rightarrow\infty}\sum^{\infty}_{k=1}a_{nk}x_{k}=\ell.
\]
Then $\ell$ is called the $\verb"A"$-limit of $\textbf{x}$. We have thus defined a method of sequential convergence, i.e. $G(\textbf{x})=\lim \verb"A"\textbf{x}$, called a matrix method or a summability matrix.

The concept of statistical convergence is a generalization of the usual notion of convergence that, for real-valued sequences, parallels the usual theory of convergence. A sequence $(\alpha _{k})$ of points in \textbf{R} is called statistically convergent to an element $\ell$ of \textbf{R}  if for each
$\varepsilon$
\[
\lim_{n\rightarrow\infty}\frac{1}{n}|\{k\leq n: |\alpha _{k}-\ell|\geq{\varepsilon}\}|=0,
\] and this is denoted by $st-\lim_{k\rightarrow\infty}\alpha _{k}=\ell$ (see \cite{FridyOnstatisticalconvergence}, \cite{MaioKocinac}, \cite{CakalliAstudyonstatisticalconvergence}, \cite{CasertaandKocinacOnstatisticalexhaustiveness}, and  \cite{CasertaMaioKocinacStatisticalConvergenceinFunctionSpaces}). This defines a method of sequential convergence, i.e. $G(\boldsymbol{\alpha}):= st-\lim_{k\rightarrow\infty}\alpha _{k}$.

 A lacunary sequence $\theta=(k_{r})$ is an increasing sequence $\theta=(k_{r})$ of positive integers such that $k_{0}=0$ and $h_{r}:k_{r}-k_{r-1}\rightarrow\infty$. The intervals determined by $\theta$ will be denoted by $I_{r}=(k_{r-1}, k_{r}]$, and the ratio $\frac{k_{r}}{k_{r-1}}$ will be abbreviated by $q_{r}$. Sums of the form $\sum^{k_{r}}_{k_{r-1}+1}|\alpha_{k}|$ frequently occur, and will often be written for convenience as $\sum^{}_{k\in{I_{r}}}|\alpha_{k}|$. Throughout this paper, we will assume that $\lim inf_{r}\; q_{r}>1$.

The notion of $N_\theta$ convergence was introduced, and studied by Freedman, Sember, and M. Raphael in \cite{FreedmanandSemberandRaphaelSomecesarotypesummabilityspaces}. Basarir and Altundag studied $\Delta$-$N_{\theta}$-asymptotically equivalent sequences in \cite{BasarirandSelmaAltundagDeltaLacunarystatistical}.

A sequence $(\alpha_{k})$ of points in \textbf{R} is called $N_\theta$-convergent to an element $\ell$ of \textbf{R} if
\[
\lim_{r\rightarrow\infty}\frac{1}{h_{r}}\sum^{}_{k\in{I_{r}}}|\alpha_{k}-\ell|=0,
\]
and it is denoted by $N_{\theta}-lim\; \alpha_{k}=\ell$. This defines a method of sequential convergence, i.e. $G(\boldsymbol{\alpha}):=N_{\theta}-lim\; \alpha_{k}$. Any convergent sequence is $N_{\theta}$-convergent, but the converse is not always true. Throughout the paper $N_{\theta}$ will denote the set of $N_\theta$ convergent sequences of points in \textbf{R}.

Using the idea of Sember and Raphael, Fridy and Orhan (\cite{FridyandOrhanlacunarystatisconvergence})  introduced lacunary statistical convergence (see also \cite{FridyandOrhanlacunarystatisticalsummability}, \cite{SavasOnlacunarystrongsigmaconvergence}, \cite{SavasandNurayOnsigmastatisticallyconvergenceanlacunarysigmastatisticallyconvergence}, \cite{SavasOnlacunarystatisticallyconvergentdoublesequencesoffuzzynumbers}, \cite{MursaleenandMohiuddineOnlacunarystatisticalconvergencewithrespecttotheintuitionisticfuzzynormedspace}, and \cite{MursaleenandAlotaibiStatisticallacunarysummabilityandaKorovkintypeapproximationtheorem}).

A real sequence $(x_{k})$ is called lacunary statistically convergent to an element $\ell$ of $\textbf{R}$ if
\[
\lim_{r\rightarrow\infty}\frac{1}{h_{r}} |\{k\in I_{r}: |x_{k}-\ell|\geq{\epsilon} \}|=0,
\]
for every $\epsilon>0$ where $I_{r}=(k_{r-1},k_{r}]$ and $k_{0}=0$,
$h_{r}:k_{r}-k_{r-1}\rightarrow \infty$ as $r\rightarrow\infty$ and
$\theta=(k_{r})$ is an increasing sequence of positive integers. For
an introduction to lacunary statistical convergence for real or
complex number sequences see \cite{FridyandOrhanlacunarystatisconvergence} and  \cite{CakalliLacunarystatisticalconvergenceintopgroups}).

For example, limit of the sequence of the ratios of Fibonacci numbers converge to the golden mean. This property ensures the regularity of lacunary sequential method obtained via the sequence of Fibonacci numbers, i.e. $\theta=(k_{r})$ is the lacunary sequence defined by writing $k_{0}=0$ and $k_{r}=F_{r+2}$ where $(F_{r})$ is the Fibonacci sequence, i.e. $F_{1}=1$, $F_{2}=1$, $F_{r}= F_{r-1} + F_{r-2}$ for $r\geq 3$.

Connor and Grosse-Erdmann \cite{ConnorandGrosseSequentialdefinitionsofcontinuityforrealfunctions} gave sequential definitions of continuity for real functions calling $G$-continuity instead of $A$-continuity and their results cover the earlier works related to $A$-continuity where a method of sequential convergence, or briefly a method, is a linear function $G$ defined on a linear subspace of $s$, denoted by $c_{G}$, into $\textbf{R}$. A sequence $\textbf{x}=(x_{n})$ is said to be $G$-convergent to $\ell$ if $\textbf{x}\in c_{G}$ and $G(\textbf{x})=\ell$. In particular, $\lim$ denotes the limit function $\lim \textbf{x}=\lim_{n}x_{n}$ on the linear space $c$ and $st-\lim$ denotes the statistical limit function $st-\lim \textbf{x}=st-\lim_{n}x_{n}$ on the linear space $st(\textbf{R})$. A function $f$ is called $G$-continuous at a point $u$ provided that whenever a sequence $\textbf{x}=(x_{n})$ of terms in the domain of $f$ is $G$-convergent to $u$, then the sequence $f(\textbf{x})=(f(x_{n}))$ is $G$-convergent to $f(u)$. A method $G$ is called regular if every convergent sequence $\textbf{x}=(x_{n})$ is $G$-convergent with $G(\textbf{x})=\lim \textbf{x}$. A method is called subsequential if whenever $\textbf{x}$ is $G$-convergent with $G(\textbf{x})=\ell$, then there is a subsequence $(x_{n_{k}})$ of $\textbf{x}$ with $\lim_{k} x_{n_{k}}=\ell$. Recently, Cakalli  gave new sequential definitions of compactness and slowly oscillating compactness in \cite{CakalliSequentialdefinitionsofcompactness}, and \cite{CakalliSlowlyoscillatingcontinuity} (see also \cite{Vallin}), respectively.

Now we recall the concepts of ward compactness, and slowly oscillating compactness: a subset $A$ of  \textbf{R} is called ward compact if whenever $(\alpha _{n})$ is a sequence of points in $A$ there is a quasi-Cauchy subsequence $\textbf{z}=(z_{k})=(\alpha_{n_{k}})$ of $(\alpha _{n})$ (\cite{CakalliForwardcontinuity}). A subset $A$ of  \textbf{R} is called slowly oscillating compact if whenever $(\alpha _{n})$ is a sequence of points in $A$ there is a slowly oscillating subsequence $\textbf{z}=(z_{k})=(\alpha_{n_{k}})$ of $(\alpha _{n})$ (\cite{CakalliSlowlyoscillatingcontinuity}).\\

 A function $f$ is called $G$-sequentially continuous at $u\in {\textbf{R}}$ if, given a sequence $\boldsymbol{\alpha}=(\alpha_{n})$ of points in \textbf{R}, $G(\boldsymbol{\alpha})=u$ implies that $G(f(\boldsymbol{\alpha}))=f(u)$.

Recently, \c{C}akall\i\; (\cite[page 594]{CakalliSequentialdefinitionsofcompactness}, \cite{CakalliOnGcontinuity}) gave a sequential definition of compactness, which is a generalization of ordinary sequential compactness, as in the following: a subset $A$ of \textbf{R} is $G$-sequentially compact if for any sequence $(\alpha _{k})$ of points in $A$ there exists a subsequence $\textbf z$ of the sequence such that $G(\textbf{z})\in{A}$. His idea enables us obtaining new kinds of compactness via most of the non-matrix sequential convergence methods, as well as all matrix sequential convergence methods.

Recently, Palladino (\cite{PalladinoOnhalfcauchysequences}) introduced a concept of upward half Cauchness, and a concept of downward half Cauchyness as in the following:

\begin{Def} \label{DefinitionHalfupwardCauchysequence} A sequence $(x_n)$ of points in $\textbf{R}$ is called upward half Cauchy if for every
$\varepsilon>0$ there exists an $n_{0}\in{\textbf{N}}$ so that $x_{n}-x_{m} <\varepsilon$ for $m \geq n \geq n_0$.
\end{Def}

Any Cauchy sequence is upward half Cauchy, so any convergent sequence is.  Any subsequence of an upward half Cauchy sequence is again upward half Cauchy.

\begin{Def} \label{DefinitionHalfdownwardCauchysequence} We say that a sequence $(x_n)$ of points in $\textbf{R}$ is downward half Cauchy if for every $\varepsilon>0$ there exists an $n_{0}\in{\textbf{N}}$ so that $x_{m}-x_{n} <\varepsilon$ for $m \geq n \geq n_0$.
\end{Def}

Any Cauchy sequence is downward half Cauchy, so any convergent sequence is. Any subsequence of a downward half Cauchy sequence is again downward half Cauchy.

In \cite{PalladinoOnhalfcauchysequences} a sequence $(x_n)$ of points in $\textbf{R}$ is called half Cauchy if the sequence is either upward half Cauchy, or downward half  Cauchy, or both and it was proved that a sequence $(x_n)$ of points in $\textbf{R}$ is Cauchy if and only if it is both upward half Cauchy and downward half Cauchy.

\maketitle

\section{Half quasi-Cauchy sequences}

First we recall that a sequence $(x_n)$ of points in $\textbf{R}$ is quasi Cauchy if $(x_{n}-x_{n+1})$ is a null sequence. Recently quasi-Cauchy sequences have been studied by Burton and Coleman \cite{BurtonColeman}, and by the author \cite{CakalliForwardcontinuity} and they obtained interesting results on uniform continuity.

Quasi-Cauchy sequences arise in diverse situations, and it is often difficult to determine whether or not they converge, and if so, to which limit. It is easy to construct a zero-one sequence such that the quasi-Cauchy average sequence does not converge. The usual constructions have a somewhat artificial feeling. Nevertheless, there are sequences which seem natural, have the quasi-Cauchy property, and do not converge. We give three interesting examples which are not sequences of averages.

\textbf{Example 1.} Let $n$ be a positive integer. In a group of $n$ people, each person selects at random and simultaneously another person of the group. All of the selected persons are then removed from the group, leaving a random number $n_{1} < n$ of people which form a new group. The new group then repeats independently the selection and removal thus described, leaving $n_{2} < n_{1}$ persons, and so forth until either one person remains, or no persons remain. Denote by $p_n$ the probability that, at the end of this iteration initiated with a group of $n$ persons, one person remains. Then the sequence
$\textbf{p} = (p_{1}, p_{2}, · · ·, p_{n},...)$  is a quasi-Cauchy sequence, and $lim p_n$ does not exist (\cite{WinklerMathematicalPuzzles}).

\textbf{Example 2.} Let $n$ be a positive integer. In a group of $n$ people, each person selects independently and at random one of three subgroups to which to
belong, resulting in three groups with random numbers $n_{1}$, $n_{2}$, $n_{3}$ of members; $n_{1} + n_{2} + n_{3} = n$. Each of the subgroups is then partitioned independently in the same manner to form three sub subgroups, and so forth. Subgroups having no members or having only one member are removed from the process. Denote by $t_{n}$ the expected value of the number of iterations up to complete removal, starting initially with a group of $n$ people. Then the sequence $(t_{1}, \frac{t_{2}}{2}, \frac{t_{3}}{3},...,\frac{t_{k}}{k},...)$ is a bounded nonconvergent quasi-Cauchy sequence.

\textbf{Example 3.} Let $\textbf{x}:= (x_{n})$ be a sequence such that for each nonnegative integer $n$, $x_{n}$ is either $0$ or $1$. For each positive integer $n$ set  $a_{n}=\frac{x_{1}+x_{2}+...+x_{n}}{n}$ . Then $a_{n}$ is the arithmetic mean average of the sequence up to time or position $n$. Clearly
for each $n$, $0 \leq a_{n} \leq 1$. $(a_{n})$ is a quasi-Cauchy sequence. That is, the sequence of averages of $0$ s and $1$ s is always a quasi-Cauchy sequence.

Using the idea of the definition of an upward half Cauchy sequence, we introduce a concept of upward half quasi-Cauchy sequence.

\begin{Def} We say that a sequence $(x_n)$ of points in $\textbf{R}$ is upward half quasi-Cauchy if for every
$\varepsilon>0$ there exists an $n_{0}\in{\textbf{N}}$ such that $x_{n}-x_{n+1} <\varepsilon$ for $n \geq n_0$.
\end{Def}

Any quasi-Cauchy sequence is upward half quasi-Cauchy, so any slowly oscillating sequence is upward half quasi-Cauchy, so any Cauchy sequence is, so any convergent sequence is. Any upward half Cauchy sequence is upward half quasi-Cauchy.

As in the case that a sequence $\textbf{x}=(x_{n})$ is Cauchy if and only if every subsequence of $\textbf{x}$ is quasi-Cauchy we have the following for upward half Cauchyness.

\begin{Thm} \label{TheohalfupwardCauchyiffeverysubsequenceisupwardquasi-Cauchy}
A sequence $\textbf{x}=(x_n)$ of points in $\textbf{R}$ is upward half Cauchy if and only if every subsequence of $\textbf{x}$ is upward half quasi-Cauchy.

\end{Thm}

\begin{pf}
It is clear that if $\textbf{x}=(x_n)$ is upward half Cauchy, then any subsequence of $\textbf{x}$ is upward half quasi-Cauchy. To prove the converse now suppose that $(x_n)$ is not upward half Cauchy so that there is a positive real number $\varepsilon_{0}$ such that for every positive integer $n$ there exist positive integers $k_{n}$ and $k_{m}$ satisfying $k_{m}>k_{n}>n$ and $x_{k_{n}}-x_{k_{m}}\geq {\varepsilon_{0}}$. Choose positive integers $k_{1}$ and $k_{2}$ satisfying $k_{2}>k_{1}>1$ and $x_{k_{1}}-x_{k_{2}}\geq {\varepsilon_{0}}$. For $n=2$ choose positive integers $k_{3}$ and $k_{4}$ satisfying $k_{4}>k_{3}>k_{2}$ and $x_{k_{3}}-x_{k_{4}}\geq {\varepsilon_{0}}$. Having inductively chosen  positive integers $k_{n}$ and $k_{n+1}$ satisfying $k_{n+1}>k_{n}>k_{n-1}$ and $x_{k_{n}}-x_{k_{n+1}}\geq {\varepsilon_{0}}$ that we can choose positive integers $k_{n+2}$ and $k_{n+3}$ satisfying $k_{n+3}>k_{n+2}>k_{n+1}$ and $x_{k_{n+2}}-x_{k_{n+3}}\geq {\varepsilon_{0}}$. Hence the subsequence $(x_{k_{n}})$ is not upward half quasi-Cauchy. Thus the proof of the theorem is completed.
\end{pf}

Now we introduce a definition of upward half Cauchy compactness of a subset of $\textbf{R}$.

\begin{Def}
A subset $E$ of $\textbf{R}$ is called upward half Cauchy compact if whenever $\textbf{x}=(x_{n})$ is a sequence of points in $E$ there is an upward half Cauchy subsequence $\textbf{z}=(z_{k})=(x_{n_{k}})$ of $\textbf{x}$.
\end{Def}

Firstly, we note that any finite subset of $\textbf{R}$ is upward half Cauchy compact, union of two upward half Cauchy  compact subsets of $\textbf{R}$ is upward half Cauchy  compact and intersection of any upward half Cauchy  compact subsets of $\textbf{R}$ is upward half Cauchy  compact. Furthermore any subset of an upward half Cauchy  compact set is upward half Cauchy compact and any bounded subset of $\textbf{R}$ is upward half Cauchy  compact. Any compact subset of $\textbf{R}$ is also upward half Cauchy  compact. We note that any slowly oscillating compact subset of $\textbf{R}$ is upward half Cauchy  compact (see [7] for the definition of slowly oscillating compactness).

\begin{Def} We say that a sequence $(x_n)$ of points in $\textbf{R}$ is downward half quasi-Cauchy if for every
$\varepsilon>0$ there exists an $n_{0}\in{\textbf{N}}$ such that $x_{n+1}-x_{n} <\varepsilon$ for $n \geq n_0$.
\end{Def}

Any quasi-Cauchy sequence is downward half quasi-Cauchy, so any slowly oscillating sequence is downward half quasi-Cauchy, so any Cauchy sequence is, so any convergent sequence is. Any downward half Cauchy sequence is downward half quasi-Cauchy.

As in the case that a sequence $\textbf{x}=(x_{n})$ is upward half Cauchy if and only if every subsequence of $\textbf{x}$ is upward half quasi-Cauchy we have the following for downward half Cauchyness.

\begin{Thm} \label{TheohalfdownwardCauchyiffeverysubsequenceisdownwardquasi-Cauchy}
A sequence $(x_n)$ of points in $\textbf{R}$ is downward half Cauchy if and only if every subsequence of $(x_n)$ is downward half quasi-Cauchy.
\end{Thm}

\begin{pf}
It is clear that if $\textbf{x}=(x_n)$ is downward half Cauchy, then any subsequence of $\textbf{x}$  is downward half quasi-Cauchy. To prove the converse now suppose that $(x_n)$ is not downward half Cauchy so that there is a positive real number $\varepsilon_{0}$ such that for every positive integer $n$ there exist positive integers $k_{n}$ and $k_{m}$ satisfying $k_{m}>k_{n}>n$ and $x_{k_{m}}-x_{k_{n}}\geq {\varepsilon_{0}}$. Choose positive integers $k_{1}$ and $k_{2}$ satisfying $k_{2}>k_{1}>1$ and $x_{k_{2}}-x_{k_{1}}\geq {\varepsilon_{0}}$. For $n=2$ choose positive integers $k_{3}$ and $k_{4}$ satisfying $k_{4}>k_{3}>k_{2}$ and $x_{k_{4}}-x_{k_{3}}\geq {\varepsilon_{0}}$. Having inductively chosen  positive integers $k_{n}$ and $k_{n+1}$ satisfying $k_{n+1}>k_{n}>k_{n-1}$ and $x_{k_{n+1}}-x_{k_{n}}\geq {\varepsilon_{0}}$ that we can choose positive integers $k_{n+2}$ and $k_{n+3}$ satisfying $k_{n+3}>k_{n+2}>k_{n+1}$ and $x_{k_{n+3}}-x_{k_{n+2}}\geq {\varepsilon_{0}}$. Thus the subsequence $(x_{k_{n}})$ is not downward half quasi-Cauchy. Hence the proof of the theorem is completed.

\end{pf}

Now we introduce a definition of downward half Cauchy compactness of a subset of $\textbf{R}$.

\begin{Def}
A subset $E$ of $\textbf{R}$ is called downward half Cauchy compact if whenever $\textbf{x}=(x_{n})$ is a sequence of points in $E$ there is a downward half Cauchy subsequence $\textbf{z}=(z_{k})=(x_{n_{k}})$ of $\textbf{x}$.
\end{Def}

Firstly, we note that any finite subset of $\textbf{R}$ is downward half Cauchy compact, union of two downward half Cauchy  compact subsets of $\textbf{R}$ is downward half Cauchy  compact and intersection of any downward half Cauchy  compact subsets of $\textbf{R}$ is downward half Cauchy compact. Furthermore any subset of an downward half Cauchy  compact set is downward half Cauchy  compact and any bounded subset of $\textbf{R}$ is downward half Cauchy  compact. Any compact subset of $\textbf{R}$ is also downward half Cauchy  compact, and the set $\textbf{N}$ is not downward half Cauchy  compact. We note that any slowly oscillating compact subset of $\textbf{R}$ is downward half Cauchy compact.

\begin{Def} We say that a sequence $(x_{n})$ of points in $\textbf{R}$ is half quasi-Cauchy if the sequence is either upward half quasi-Cauchy, or downward half quasi-Cauchy, or both.
\end{Def}

As we expected, we have the following:

\begin{Thm}  A sequence $(x_{n})$ of points in $\textbf{R}$ is quasi-Cauchy if and only if it is both upward half quasi-Cauchy and downward half quasi-Cauchy.
\end{Thm}
\begin{pf}
Suppose that a sequence $(x_n)$ of points in $\textbf{R}$ is quasi-Cauchy and choose an arbitrary $\varepsilon> 0$. Since $(x_n)$ is quasi-Cauchy there exists a positive integer $n_{0}$ such that $|x_{n}-x_{n+1}|<\varepsilon$ for $n>n_{0}$. Now we obtain that  $x_{n}-x_{n+1}\leq |x_{n}-x_{n+1}|<\varepsilon$ and $x_{n+1}-x_{n}\leq |x_{n}-x_{n+1}|<\varepsilon$ for $n>n_{0}$. Thus $(x_n)$ is both upward half quasi-Cauchy and downward
half quasi-Cauchy. Now suppose $(x_{n})$ is both upward half quasi-Cauchy and downward half quasi-Cauchy. Choose an arbitrary $\varepsilon> 0$. Since $(x_{n})$ is upward half quasi-Cauchy, there exists a positive integer $n_{1}$ such that $x_{n}-x_{n+1} <\varepsilon$ for $n \geq n_{1}$. Since $(x_n)$ is downward half quasi-Cauchy, there exists a positive integer $n_{2}$ such that $x_{n+1}-x_{n} <\varepsilon$ for $n \geq n_{2}$. Write $n_{0}=max \{n_{1}, n_{2}\}$. Then $-\varepsilon < x_{n}-x_{n+1}  <\varepsilon$ for $n>n_{0}$, so  $|x_{n}-x_{n+1}|<\varepsilon$ for $n>n_{0}$. This completes the proof of the theorem.
\end{pf}

A real function $f$ is continuous if and only if, for each point $\ell$ in the domain, $\lim_{n\rightarrow\infty}f(x_{n})=f(\ell)$ whenever $\lim_{n\rightarrow\infty}x_{n}=\ell$. This is equivalent to the statement that $(f(x_{n}))$ is a convergent sequence whenever $(x_{n})$ is. This is also equivalent to the statement that $(f(x_{n}))$ is a Cauchy sequence whenever $(x_{n})$ is Cauchy. These well known results for continuity for real functions in terms of sequences might suggest to us giving a new type continuity, namely, upward half Cauchy continuity:

\begin{Def}
A function $f$ is called upward half Cauchy continuous on $E$ if the sequence $(f(x_{n}))$ is upward half Cauchy whenever $\textbf{x}=(x_{n})$ is an upward half Cauchy sequence of points in $E$.
\end{Def}

We note that sum of two upward half Cauchy  continuous functions is upward  half Cauchy continuous and composite of two upward half Cauchy  continuous functions is upward  half Cauchy  continuous.

In connection with upward half Cauchy sequences and convergent sequences the problem arises to investigate the following types of  "continuity" of functions on $\textbf{R}$ where $C^{+}$ denotes the set of all upward half Cauchy sequences of points in $\textbf{R}$.

\begin{description}
\item[($C^{+} $)] $(x_{n}) \in {C^{+}} \Rightarrow (f(x_{n})) \in {C^{+}}$
\item[($C^{+} c$)] $(x_{n}) \in {C^{+}} \Rightarrow (f(x_{n})) \in {c}$
\item[$(c)$] $(x_{n}) \in {c} \Rightarrow (f(x_{n})) \in {c}$
\item[$(cC^{+})$] $(x_{n}) \in {c} \Rightarrow (f(x_{n})) \in {C^{+}}$
\end{description}

We see that $(C^{+})$ is upward half Cauchy  continuity of $f$ and $(c)$ states the ordinary continuity of $f$. It is easy to see that $(C^{+} c)$ implies $(C^{+})$, and $(C^{+} )$ does not imply $(C^{+} c)$;  and $(C^{+})$ implies $(cC^{+})$, and $c(C^{+})$ does not imply $(C^{+})$; $(C^{+} c)$ implies $(c)$ and $(c)$ does not imply $(C^{+} c)$; and $(c)$ is equivalent to $(cC^{+})$.

We see that (c) can be replaced by statistical continuity, i.e.  $st-\lim_{n\rightarrow \infty} f(x_{n})=f(\ell)$  whenever $\textbf{x}=(x_{n})$ is a statistically convergent sequence with $st-\lim_{n\rightarrow \infty} x_{n}=\ell$.

Now we give the implication $(C^{+})$ implies $(c)$, i.e. any upward half Cauchy continuous function is continuous in the ordinary sense.

\begin{Thm} If $f$ is upward half Cauchy continuous on a subset $E$ of $\textbf{R}$, then it is continuous on $E$ in the ordinary sense.
\end{Thm}
\begin{pf}

Let $(x_{n})$ be any convergent sequence with $\lim_{k\rightarrow\infty}x_{k}=\ell$. Then
$$(x_{1}, \ell , x_{2}, \ell,..., x_{n}, \ell,...)$$ is also convergent to $\ell$. So it is upward half Cauchy. Hence
$$(f(x_{1}),f(\ell),f(x_{2}),f(\ell),...,f(x_{n}),f(\ell),...)$$ is upward half Cauchy, therefore it is upward half quasi-Cauchy. Let $\varepsilon>0$. There exists an $n_{0}\in{\textbf{N}}$ so that $f(x_{n})-\ell <\varepsilon$ and $\ell-f(x_{n})<\varepsilon$ for $n \geq n_0$. Hence $|f(x_{n})-\ell| <\varepsilon$ for $n \geq n_0$. It follows from this that the sequence $(f(x_{n}))$ converges to $f(\ell)$. This completes the proof of the theorem.
\end{pf}

Now we state the following result related to statistical continuity and upward half Cauchy continuity.

\begin{Cor} If $f$ is upward half Cauchy  continuous, then it is statistically continuous.
\end{Cor}

\begin{Thm} \label{ThehalfquasiCauchycontinousimageofhalfquasi-Cauchycompactsubsetis}
 Upward half Cauchy  continuous image of any upward half Cauchy compact subset of $\textbf{R}$ is upward half Cauchy  compact.
\end{Thm}
\begin{pf}
Write $y_{n}=f(x_{n})$ where $x_{n}\in {E}$ for each $n \in{\textbf{N}}$. Upward half Cauchy  compactness of $E$ implies that there is an upward half Cauchy subsequence $\textbf{z}=(z_{k})=(x_{n_{k}})$ of $\textbf{x}$. Write $(t_{k})=f(\textbf{z})=(f(z_{k}))$. $(t_{k})$ is an upward half Cauchy subsequence of the sequence $f(\textbf{x})$. This completes the proof of the theorem.
\end{pf}

\begin{Cor} Upward half Cauchy continuous image of any compact subset of $\textbf{R}$ is compact.
\end{Cor}

It is a well known result that uniform limit of a sequence of continuous functions is continuous. This is also true in case upward half Cauchy  continuity, i.e. uniform limit of a sequence of upward half Cauchy  continuous functions is upward half Cauchy  continuous.

\begin{Thm} If $(f_{n})$ is a sequence of upward half Cauchy  continuous functions defined on a subset $E$ of $\textbf{R}$ and $(f_{n})$ is uniformly convergent to a function $f$, then $f$ is upward half Cauchy  continuous on $E$.
\end{Thm}
\begin{pf}  Let $\varepsilon > 0$. There exists a positive integer $N$ such that $|f_{n}(x)-f(x)|<\frac{\varepsilon}{3}$ for all $x \in {E}$ whenever $n\geq N$. Take any upward half Cauchy sequence $(x_{n})$ of points in $E$. Since $f_{N}$ is upward half Cauchy continuous, there exists a positive integer $N_{1}$, depending on $\varepsilon$, and greater than $N$ such that $f_{N}(x_{n})-f_{N}(x_{m})<\frac{\varepsilon}{3}$ for $m \geq n \geq N_{1}$. Now for $m \geq n \geq N_{1}$  we have $$ f(x_{n})-f(x_{m})=f(x_{n})-f_{N}(x_{n})+f_{N}(x_{n})-f_{N}(x_{m})+f_{N}(x_{m})-f(x_{m}) $$ $$\leq
f(x_{n})-f_{N}(x_{n})+\frac{\varepsilon}{3} +f_{N}(x_{m})-f(x_{m})$$
$$\leq |f(x_{n})-f_{N}(x_{n})|+ \frac{\varepsilon}{3}+|f_{N}(x_{m})-f(x_{m})|\leq \frac{\varepsilon}{3} +\frac{\varepsilon}{3}+ \frac{\varepsilon}{3}= \varepsilon.$$
\end{pf}

Now we give the definition of downward half Cauchy continuity of a real function analogous to the definition of upward half Cauchy continuity.

\begin{Def}
A function $f$ is called downward half Cauchy continuous on $E$ if the sequence $(f(x_{n}))$ is downward half Cauchy whenever $\textbf{x}=(x_{n})$ is a downward half Cauchy sequence of points in $E$.
\end{Def}

We note that sum of two downward half Cauchy  continuous functions is downward  half Cauchy continuous and composite of two downward half Cauchy continuous functions is downward  half Cauchy  continuous.

In connection with downward half Cauchy sequences and convergent sequences the problem arises to investigate the following types of  "continuity" of functions on $\textbf{R}$ where $C^{-}$ denotes the set of all downward half Cauchy sequences of points in $\textbf{R}$.

\begin{description}
\item[($C^{-} $)] $(x_{n}) \in {C^{-}} \Rightarrow (f(x_{n})) \in {C^{-}}$
\item[($C^{-} c$)] $(x_{n}) \in {C^{-}} \Rightarrow (f(x_{n})) \in {c}$
\item[$(c)$] $(x_{n}) \in {c} \Rightarrow (f(x_{n})) \in {c}$
\item[$(c C^{-})$] $(x_{n}) \in {c} \Rightarrow (f(x_{n})) \in {C^{-}}$
\end{description}

We see that $(C^{-})$ is downward half Cauchy  continuity of $f$. It is easy to see that $(C^{-} c)$ implies $(C^{-})$, and $(C^{-} )$ does not imply $(C^{-} c)$;  and $(C^{-})$ implies $(c C^{-})$, and $(c C^{-})$ does not imply $(C^{-})$; $(C^{-} c)$ implies $(c)$ and $(c)$ does not imply $(C^{-}c)$; and $(c)$ is equivalent to $(c C^{-})$.

Now we give the implication $(C^{-})$ implies $(c)$, i.e. any upward half Cauchy continuous function is continuous in the ordinary sense.

\begin{Thm} If $f$ is downward half Cauchy continuous on a subset $E$ of $\textbf{R}$, then it is continuous on $E$ in the ordinary sense.
\end{Thm}
\begin{pf}

Let $(x_{n})$ be any convergent sequence with $\lim_{k\rightarrow\infty}x_{k}=\ell$. Then
$$(x_{1}, \ell , x_{2}, \ell,..., x_{n}, \ell,...)$$ is also convergent to $\ell$. So it is downward half Cauchy. Hence
$$(f(x_{1}),f(\ell),f(x_{2}),f(\ell),...,f(x_{n}),f(\ell),...)$$ is downward half Cauchy, therefore it is downward half quasi-Cauchy. Let $\varepsilon>0$. There exists an $n_{0}\in{\textbf{N}}$ so that $f(x_{n})-\ell <\varepsilon$ and $\ell-f(x_{n})<\varepsilon$ for $n \geq n_0$. Hence $|f(x_{n})-\ell| <\varepsilon$ for $n \geq n_0$. It follows from this that the sequence $(f(x_{n}))$ converges to $f(\ell)$. This completes the proof of the theorem.
\end{pf}

Now we state the following result related to statistical continuity and downward half Cauchy continuity.

\begin{Cor} If $f$ is downward half Cauchy  continuous, then it is statistically continuous.
\end{Cor}

\begin{Thm} \label{ThedownhalfquasiCauchycontinousimageofhalfquasi-Cauchycompactsubsetis}
Downward half Cauchy  continuous image of any downward half Cauchy compact subset of $\textbf{R}$ is downward half Cauchy  compact.
\end{Thm}
\begin{pf}
Write $y_{n}=f(x_{n})$ where $x_{n}\in {E}$ for each $n \in{\textbf{N}}$. Downward half Cauchy  compactness of $E$ implies that there is a downward half Cauchy subsequence $\textbf{z}=(z_{k})=(x_{n_{k}})$ of $\textbf{x}$. Write $(t_{k})=f(\textbf{z})=(f(z_{k}))$. $(t_{k})$ is a downward half Cauchy subsequence of the sequence $f(\textbf{x})$. This completes the proof of the theorem.
\end{pf}

\begin{Cor} Downward half Cauchy continuous image of any compact subset of $\textbf{R}$ is compact.
\end{Cor}

It is a well known result that uniform limit of a sequence of continuous functions is continuous. This is also true in the case of downward half Cauchy  continuity, i.e. uniform limit of a sequence of downward half Cauchy  continuous functions is downward half Cauchy  continuous.

\begin{Thm} If $(f_{n})$ is a sequence of downward half Cauchy  continuous functions defined on a subset $E$ of $\textbf{R}$ and $(f_{n})$ is uniformly convergent to a function $f$, then $f$ is downward half Cauchy  continuous on $E$.
\end{Thm}
\begin{pf}  Let $\varepsilon > 0$. Then there exists a positive integer $N$ such that $|f_{n}(x)-f(x)|<\frac{\varepsilon}{3}$ for all $x \in {E}$ whenever $n\geq N$. Take any downward half Cauchy sequence $(x_{n})$ of points in $E$. As $f_{N}$ is downward half Cauchy continuous, the sequence $(f_{N}(x_{n}))$ is a downward half Cauchy sequence, so there exists a positive integer $N_{1}$, depending on $\varepsilon$, and greater than $N$ such that $f_{N}(x_{m})-f_{N}(x_{n})<\frac{\varepsilon}{3}$ for $m \geq  n\geq N_{1}$. Now for $m \geq n\geq N_{1}$  we have $$ f(x_{m})-f(x_{n})=f(x_{m})-f_{N}(x_{m})+f_{N}(x_{m})-f_{N}(x_{n})+f_{N}(x_{n})-f(x_{n}) $$ $$\leq
f(x_{m})-f_{N}(x_{m})+\frac{\varepsilon}{3} +f_{N}(x_{n})-f(x_{n})$$
$$\leq |f(x_{m})-f_{N}(x_{m})|+ \frac{\varepsilon}{3}+|f_{N}(x_{n})-f(x_{n})|\leq \frac{\varepsilon}{3} +\frac{\varepsilon}{3}+ \frac{\varepsilon}{3}= \varepsilon.$$
\end{pf}

Now we introduce a definition of upward half quasi-Cauchy compactness of a subset of $\textbf{R}$.

\begin{Def}
A subset $E$ of $\textbf{R}$ is called upward half quasi-Cauchy compact if whenever $\textbf{x}=(x_{n})$ is a sequence of points in $E$ there is an upward half quasi-Cauchy subsequence $\textbf{z}=(z_{k})=(x_{n_{k}})$ of $\textbf{x}$.
\end{Def}

Firstly, we note that any finite subset of $\textbf{R}$ is upward half quasi-Cauchy compact, union of two upward half quasi-Cauchy  compact subsets of $\textbf{R}$ is upward half quasi-Cauchy  compact and intersection of any upward half quasi-Cauchy  compact subsets of $\textbf{R}$ is upward half quasi-Cauchy  compact. Furthermore any subset of an upward half quasi-Cauchy  compact set is upward half quasi-Cauchy compact and any bounded subset of $\textbf{R}$ is upward half quasi-Cauchy  compact. Any compact subset of $\textbf{R}$ is also upward half quasi-Cauchy  compact, and the set of negative integers is not upward half quasi-Cauchy  compact. We note that any slowly oscillating compact subset of $\textbf{R}$ is upward half quasi-Cauchy  compact.

Now we give the concept of upward half quasi-Cauchy continuity of a function defined on a subset of $\textbf{R}$ to $\textbf{R}$.

\begin{Def}
A function $f$ is called upward half quasi-Cauchy continuous on $E$ if the sequence $(f(x_{n}))$ is upward half quasi-Cauchy whenever $\textbf{x}=(x_{n})$ is an upward half quasi-Cauchy sequence of points in $E$.
\end{Def}

We see that sum of two upward half quasi-Cauchy  continuous functions is upward  half quasi-Cauchy continuous and composite of two upward half quasi-Cauchy  continuous functions is upward  half quasi-Cauchy  continuous.

In connection with upward half quasi-Cauchy sequences and convergent sequences the problem arises to investigate the following types of  "continuity" of functions on $\textbf{R}$.

\begin{description}
\item[($\delta^{+} $)] $(x_{n}) \in {\Delta^{+}} \Rightarrow (f(x_{n})) \in {\Delta^{+}}$
\item[($\delta^{+} c$)] $(x_{n}) \in {\Delta^{+}} \Rightarrow (f(x_{n})) \in {c}$
\item[$(c)$] $(x_{n}) \in {c} \Rightarrow (f(x_{n})) \in {c}$
\item[$(c\delta^{+})$] $(x_{n}) \in {c} \Rightarrow (f(x_{n})) \in {\Delta^{+}}$
\end{description}

We see that $(\delta^{+})$ is upward half quasi-Cauchy  continuity of $f$. It is easy to see that $(\delta^{+} c)$ implies $(\delta^{+})$, and $(\delta^{+})$ does not imply $(\delta^{+} c)$;  and $(\delta^{+})$ implies $(c\delta^{+})$, and $(c\delta^{+})$ does not imply $(\delta^{+})$; $(\delta^{+} c)$ implies $(c)$ and $(c)$ does not imply $(\delta^{+} c)$; and $(c)$ is equivalent to $(c\delta^{+})$.

Now we prove that ($\delta^{+} $) implies $(c)$ in the following:

\begin{Thm} If $f$ is upward half quasi-Cauchy continuous on a subset $E$ of $\textbf{R}$, then it is continuous on $E$ in the ordinary sense.
\end{Thm}
\begin{pf}

Let $(x_{n})$ be any convergent sequence with $\lim_{k\rightarrow\infty}x_{k}=\ell$. Then
$$(x_{1}, \ell , x_{2}, \ell,..., x_{n}, \ell,...)$$ is also convergent to $\ell$. So it is upward half quasi-Cauchy. Hence
$$(f(x_{1}),f(\ell),f(x_{2}),f(\ell),...,f(x_{n}),f(\ell),...)$$ is upward half quasi-Cauchy. Let $\varepsilon>0$. There exists an $n_{0}\in{\textbf{N}}$ so that $f(x_{n})-\ell <\varepsilon$ and $\ell-f(x_{n})<\varepsilon$ for $n \geq n_0$. Hence $|f(x_{n})-\ell| <\varepsilon$ for $n \geq n_0$. It follows from this that the sequence $(f(x_{n}))$ converges to $f(\ell)$. This completes the proof of the theorem.
\end{pf}

Now we state the following result related to statistical continuity and upward half quasi-Cauchy continuity.

\begin{Cor} If $f$ is upward half quasi-Cauchy  continuous, then it is statistically continuous.
\end{Cor}

\begin{Thm} \label{ThehalfquasiCauchycontinousimageofhalfquasi-Cauchycompactsubsetis}
 Upward half quasi-Cauchy continuous image of any upward half quasi-Cauchy compact subset of $\textbf{R}$ is upward half quasi-Cauchy  compact.
\end{Thm}
\begin{pf}
Write $y_{n}=f(x_{n})$ where $x_{n}\in {E}$ for each $n \in{\textbf{N}}$. upward half quasi-Cauchy  compactness of $E$ implies that there is an upward half quasi-Cauchy subsequence of the sequence of $\textbf{x}$. Write $(t_{k})=f(\textbf{z})=(f(z_{k}))$. $(t_{k})$ is an upward half quasi-Cauchy subsequence of the sequence $f(\textbf{x})$. This completes the proof of the theorem.
\end{pf}

\begin{Cor} Upward half quasi-Cauchy continuous image of any compact subset of $\textbf{R}$ is compact.
\end{Cor}

It is a well known result that uniform limit of a sequence of continuous functions is continuous. This is also true in case upward half quasi-Cauchy  continuity, i.e. uniform limit of a sequence of upward half quasi-Cauchy  continuous functions is upward half quasi-Cauchy  continuous.

\begin{Thm} If $(f_{n})$ is a sequence of upward half quasi-Cauchy continuous functions defined on a subset $E$ of $\textbf{R}$ and $(f_{n})$ is uniformly convergent to a function $f$, then $f$ is upward half quasi-Cauchy continuous on $E$.
\end{Thm}
\begin{pf}  Let $\varepsilon > 0$. Then there exists a positive integer $N$ such that $|f_{n}(x)-f(x)|<\frac{\varepsilon}{3}$ for all $x \in {E}$ whenever $n\geq N$. Take any upward half quasi-Cauchy sequence $(x_{n})$ of points in $E$. Since $f_{N}$ is an upward half quasi-Cauchy continuous function, there exists a positive integer $N_{1}$, depending on $\varepsilon$, and greater than $N$ such that $$f_{N}(x_{n})-f_{N}(x_{n+1})<\frac{\varepsilon}{3}$$ for $n\geq N_{1}$. Now for $n\geq N_{1}$  we have $$ f(x_{n})-f(x_{n+1})=f(x_{n})-f_{N}(x_{n})+f_{N}(x_{n})-f_{N}(x_{n+1})+f_{N}(x_{n+1})-f(x_{n+1}) $$ $$\leq
f(x_{n})-f_{N}(x_{n})+\frac{\varepsilon}{3} +f_{N}(x_{n+1})-f(x_{n+1})$$
$$\leq |f(x_{n})-f_{N}(x_{n})|+ \frac{\varepsilon}{3}+|f_{N}(x_{n+1})-f(x_{n+1})|\leq \frac{\varepsilon}{3} +\frac{\varepsilon}{3}+ \frac{\varepsilon}{3}= \varepsilon.$$
\end{pf}

Now we introduce a definition of downward half quasi-Cauchy compactness of a subset of $\textbf{R}$.

\begin{Def}
A subset $E$ of $\textbf{R}$ is called downward half quasi-Cauchy compact if whenever $\textbf{x}=(x_{n})$ is a sequence of points in $E$ there is an downward half quasi-Cauchy subsequence $\textbf{z}=(z_{k})=(x_{n_{k}})$ of $\textbf{x}$.
\end{Def}

Firstly, we note that any finite subset of $\textbf{R}$ is downward half quasi-Cauchy compact, union of two downward half quasi-Cauchy compact subsets of $\textbf{R}$ is downward half quasi-Cauchy compact and intersection of any downward half quasi-Cauchy  compact subsets of $\textbf{R}$ is downward half quasi-Cauchy  compact. Furthermore any subset of an downward half quasi-Cauchy  compact set is downward half quasi-Cauchy compact and any bounded subset of $\textbf{R}$ is upward half quasi-Cauchy  compact. Any compact subset of $\textbf{R}$ is also upward half quasi-Cauchy  compact, and the set of negative integers is not downward half quasi-Cauchy compact. We note that any slowly oscillating compact subset of $\textbf{R}$ is downward half quasi-Cauchy compact.

Now we give the definition of downward half quasi-Cauchy continuity of a real function in the following:

\begin{Def}
A function $f$ is called downward half Cauchy continuous on $E$ if the sequence $(f(x_{n}))$ is downward half Cauchy whenever $\textbf{x}=(x_{n})$ is an downward half Cauchy sequence of points in $E$.
\end{Def}

We note that sum of two downward half quasi-Cauchy continuous functions is downward half quasi-Cauchy continuous and composite of two downward half quasi-Cauchy  continuous functions is downward half quasi-Cauchy  continuous.

In connection with downward half quasi-Cauchy sequences and convergent sequences the problem arises to investigate the following types of  "continuity" of functions on $\textbf{R}$.

\begin{description}
\item[($\delta^{-} $)] $(x_{n}) \in {\Delta^{-}} \Rightarrow (f(x_{n})) \in {\Delta^{-}}$
\item[($\delta^{-} c$)] $(x_{n}) \in {\Delta^{-}} \Rightarrow (f(x_{n})) \in {c}$
\item[$(c)$] $(x_{n}) \in {c} \Rightarrow (f(x_{n})) \in {c}$
\item[$(c\delta^{-})$] $(x_{n}) \in {c} \Rightarrow (f(x_{n})) \in {\Delta^{-}}$
\end{description}

We see that $(\delta^{-})$ is downward half quasi-Cauchy  continuity of $f$ and $(c)$ again states the ordinary continuity of $f$. It is easy to see that $(\delta^{-} c)$ implies $(\delta^{-})$, and $(\delta^{-})$ does not imply $(\delta^{-} c)$;  and $(\delta^{-})$ implies $(c\delta^{-})$, and $(c\delta^{-})$ does not imply $(\delta^{-})$; $(\delta^{-} c)$ implies $(c)$ and $(c)$ does not imply $(\delta^{-} c)$; and $(c)$ is equivalent to $(c\delta^{-})$.

Now we give the implication $(\delta^{-})$ implies $(c)$, i.e. any downward half quasi-Cauchy continuous function is continuous in the ordinary sense.

\begin{Thm} If $f$ is downward half quasi-Cauchy continuous on a subset $E$ of $\textbf{R}$, then it is continuous on $E$ in the ordinary sense.
\end{Thm}
\begin{pf}

Let $(x_{n})$ be any convergent sequence with $\lim_{k\rightarrow\infty}x_{k}=\ell$. Then
$$(x_{1}, \ell , x_{2}, \ell,..., x_{n}, \ell,...)$$ is also convergent to $\ell$. So it is downward half quasi-Cauchy. Hence
$$(f(x_{1}),f(\ell),f(x_{2}),f(\ell),...,f(x_{n}),f(\ell),...)$$ is downward half quasi-Cauchy. Let $\varepsilon>0$. There exists an $n_{0}\in{\textbf{N}}$ so that $f(x_{n})-\ell <\varepsilon$ and $\ell-f(x_{n})<\varepsilon$ for $n \geq n_0$. Hence $|f(x_{n})-\ell| <\varepsilon$ for $n \geq n_0$. It follows from this that the sequence $(f(x_{n}))$ converges to $f(\ell)$. This completes the proof of the theorem.
\end{pf}

Now we state the following result related to statistical continuity and downward half quasi-Cauchy continuity.

\begin{Cor} If $f$ is downward half quasi-Cauchy continuous, then it is statistically continuous.
\end{Cor}

\begin{Thm} \label{ThehalfquasiCauchycontinousimageofhalfquasi-Cauchycompactsubsetis}
 Downward half quasi-Cauchy  continuous image of any downward half quasi-Cauchy compact subset of $\textbf{R}$ is downward half quasi-Cauchy  compact.
\end{Thm}
\begin{pf}
Write $y_{n}=f(x_{n})$ where $x_{n}\in {E}$ for each $n \in{\textbf{N}}$. downward half quasi-Cauchy  compactness of $E$ implies that there is a downward half quasi-Cauchy subsequence $\textbf{z}=(z_{k})=(x_{n_{k}})$ of $\textbf{x}$. Write $(t_{k})=f(\textbf{z})=(f(z_{k}))$. Then $(t_{k})$ is a downward half quasi-Cauchy subsequence of the sequence $f(\textbf{x})$. This completes the proof of the theorem.
\end{pf}

\begin{Cor} Downward half quasi-Cauchy continuous image of any compact subset of $\textbf{R}$ is compact.
\end{Cor}

It is a well known result that uniform limit of a sequence of continuous functions is continuous. This is also true in case downward half quasi-Cauchy  continuity, i.e. uniform limit of a sequence of downward half quasi-Cauchy  continuous functions is downward half quasi-Cauchy  continuous.

\begin{Thm} If $(f_{n})$ is a sequence of downward half quasi-Cauchy  continuous functions defined on a subset $E$ of $\textbf{R}$ and $(f_{n})$ is uniformly convergent to a function $f$, then $f$ is downward half quasi-Cauchy  continuous on $E$.
\end{Thm}
\begin{pf}  Let $\varepsilon > 0$. Then there exists a positive integer $N$ such that $|f_{n}(x)-f(x)|<\frac{\varepsilon}{3}$ for all $x \in {E}$  whenever $n\geq N$. Take any downward half quasi-Cauchy sequence $(x_{n})$ of points in $E$. Since $f_{N}$ is downward half quasi-Cauchy continuous, there exists a positive integer $N_{1}$, depending on $\varepsilon$, and greater than $N$ such that $$f_{N}(x_{n+1})-f_{N}(x_{n})<\frac{\varepsilon}{3}$$ for $n\geq N_{1}$. Now for $n\geq N_{1}$  we have $$ f(x_{n+1})-f(x_{n})=f(x_{n+1})-f_{N}(x_{n+1})+f_{N}(x_{n+1})-f_{N}(x_{n})+f_{N}(x_{n})-f(x_{n}) $$ $$\leq
f(x_{n+1})-f_{N}(x_{n+1})+\frac{\varepsilon}{3} +f_{N}(x_{n})-f(x_{n})$$
$$\leq |f(x_{n+1})-f_{N}(x_{n+1})|+ \frac{\varepsilon}{3}+|f_{N}(x_{n})-f(x_{n})|\leq \frac{\varepsilon}{3} +\frac{\varepsilon}{3}+ \frac{\varepsilon}{3}= \varepsilon.$$
\end{pf}

For a further study, we suggest to investigate upward and downward half quasi-Cauchy sequences of fuzzy points, and  upward and downward half continuity for the fuzzy functions (see \cite{CakalliandPratul} for the definitions and  related concepts in fuzzy setting). However due to the change in settings, the definitions and methods of proofs will not always be analogous to those of the present work.

\end{document}